\documentstyle[11pt,fullpage,twoside,amssym,saks]{article}

\def\Xint#1{\mathchoice
   {\XXint\displaystyle\textstyle{#1}}%
   {\XXint\textstyle\scriptstyle{#1}}%
   {\XXint\scriptstyle\scriptscriptstyle{#1}}%
   {\XXint\scriptscriptstyle\scriptscriptstyle{#1}}%
   \!\int}
\def\XXint#1#2#3{{\setbox0=\hbox{$#1{#2#3}{\int}$}
     \vcenter{\hbox{$#2#3$}}\kern-.5\wd0}}
\def\dashint{\Xint-}

\begin{document}
\markboth{\centerline{B. Bojarski, V. Gutlyanski\v\i{} and V.
Ryazanov}} {\centerline{ON DIRICHLET PROBLEM FOR BELTRAMI
EQUATIONS}}
\def\kohta #1 #2\par{\par\noindent\rlap{#1)}\hskip30pt
\hangindent30pt #2\par}
\def\esssup{\operatornamewithlimits{ess\,sup}}
\def\tomes{\mathop{\longrightarrow}\limits^{mes}}
\def\ts{\textstyle}
\def\I{\roman{Im}}
\def\mes{\mbox{\rm mes}}
\def\Rm{{{\Bbb R}^m}}
\def\Rn{{{\Bbb R}^n}}
\def\Rk{{{\Bbb R}^k}}
\def\R3{{{\Bbb R}^3}}
\def\lR{{\overline {{\Bbb R}}}}
\def\lC{{\overline {{\Bbb C}}}}
\def\lQ{{\overline {{Q}}}}
\def\lz{{\overline {{z}}}}
\def\lD{{\overline {{D}}}}
\def\lR{{\overline {{\Bbb R}}}}
\def\lRn{{\overline {{\Bbb R}^n}}}
\def\lRm{{\overline {{\Bbb R}^m}}}
\def\lBn{{\overline {{\Bbb B}^n}}}
\def\Bn{{{\Bbb B}^n}}
\def\R{{\Bbb R}}
\def\N{{\Bbb N}}
\def\Z{{\Bbb Z}}
\def\C{{\Bbb C}}
\def\B{{\Bbb B}}
\def\Di{{\Bbb D}}
\def\e{{\varepsilon}}
\def\L{{\Lambda}}
\def\l{{\lambda}}
\def\f{{\varphi}}
\def\F{{\Phi}}
\def\x{{\chi}}
\def\d{{\delta }}
\def\D{{\Delta }}
\def\c{{\circ }}
\def\tg{{\tilde{\gamma}}}
\def\a{{\alpha }}
\def\b{{\beta }}
\def\p{{\psi }}
\def\m{{\mu }}
\def\n{{\nu }}
\def\r{{\rho }}
\def\t{{\tau }}
\def\S{{\Sigma }}
\def\O{{\Omega }}
\def\o{{\omega }}
\def\s{{\sigma }}

\def\z{{\zeta }}
\def\L{{{Log}}}
\def\E{{{Exp}}}
\def\g{{\gamma }}
\def\G{{\Gamma }}
\def\D{{\Delta }}
\let\text=\mbox
\let\Cal=\cal

\def\cc{\setcounter{equation}{0}
\setcounter{figure}{0}\setcounter{table}{0}}

\overfullrule=0pt

\def\eb{\begin{eqnarray}}
\def\ee{\end{eqnarray}}
\def\ebnn{\begin{eqnarray*}}
\def\eenn{\end{eqnarray*}}
\def\db{\begin{displaystyle}}
\def\de{\end{displaystyle}}
\def\tb{\begin{textstyle}}
\def\te{\end{textstyle}}
\def\exb{\begin{ex}}
\def\exe{\end{ex}}
\def\bth{\begin{theo}}
\def\eth{\end{theo}}
\def\bcor{\begin{corol}}
\def\ecor{\end{corol}}
\def\blem{\begin{lemma}}
\def\elem{\end{lemma}}
\def\brem{\begin{rem}}
\def\erem{\end{rem}}
\def\bpr{\begin{propo}}
\def\epr{\end{propo}}
\title{{\bf ON DIRICHLET PROBLEM FOR BELTRAMI\\ EQUATIONS WITH TWO CHARACTERISTICS}}

\author{{\bf B. Bojarski, V. Gutlyanski\v\i{} and V. Ryazanov}}
\date{\today \hskip 4mm ({\tt BGR-DIRICHLET-ARXIV-02-11-2012.tex})}
\maketitle

\large \abstract We establish a series of criteria on the existence
of regular solutions for the Dirichlet problem to general degenerate
Beltrami equations ${\overline {\partial}}f\, =\, \mu {\partial
f}+\nu {\overline {\partial f} }$ in arbitrary Jordan domains in
$\C$.
\endabstract


\large
\cc
\section{Introduction} Let $D$ be a domain in the complex plane
$\C.$ Throughout this paper we use the notations $z=x+iy,$ $B(z_0,r)\ \colon
=\{ z\in\C: |z-z_0|<r\}$ for $z_0\in\C$ and $r>0$, $\B(r)\ \colon
=B(0,r)$, $\B\ \colon =\B(1)$, and $\overline{\C}\ \colon
=\C\cup{\infty}.$
\medskip
\par
\noindent
The purpose of this paper is to study the Dirichlet problem
\begin{equation}\label{Dirichlet}
\left\{\begin{array}{ccc}
f_{\overline{z}}\, =\, \mu
(z)\cdot{f_z} + \nu (z)\cdot \overline {f_z},\,\,\, &z\in D, \\
\lim\limits_{z\to\zeta}{\rm Re}\,f(z)=\varphi(\zeta), &\forall\
\zeta\in\partial D,
\end{array}\right.
\end{equation}
in a Jordan domain $D$ of the complex plane ${\Bbb C}$ with
continuous boundary data $\varphi(\zeta)\not\equiv{\rm const}.$
Here $\mu(z)$ and $\nu(z)$ stand for measurable coefficients
satisfying the inequality $|\mu (z)|+ |\nu (z)| < 1$ a.e. in $D.$
The degeneracy of the ellipticity for the Beltrami equations
\begin{equation}\label{Beltrami}
f_{\overline{z}}\, =\, \mu
(z)\cdot{f_z} + \nu (z)\cdot \overline {f_z}
\end{equation}
 is controlled by the dilatation
coefficient
\begin{equation} \label{eq1.4}
K_{\mu , \nu}(z)\ \colon =\ \frac{1+|\mu (z)|+|\nu (z)|}{1-|\mu
(z)|-|\nu (z)|}\ \in L^1_{\rm loc}.
\end{equation}
We will look for a solution as a continuous, discrete and open mapping $f:D\to{\Bbb
C}$ of the Sobolev class $W_{\rm loc}^{1,1}$ and such that the Jacobian $J_f(z)\neq0$ a.e. in $D.$
Such a solution we will call a {\bf regular
solution}  of the Dirichlet problem (\ref{Dirichlet})
 in a domain $D.$
 \par
 Recall that a
mapping $f:D\to{\Bbb C}$ is called {\bf discrete} if the preimage
$f^{-1}(y)$ consists of isolated points for every $y\in{\Bbb C}$,
and {\bf open} if $f$ maps every open set $U\subseteq D$ onto an
open set in ${\Bbb C}$.

For the uniformly elliptic case, i.e. when $K_{\mu,\nu}(z)\leq K<\infty$ a.e. in $D$
the
Dirichlet problem was  studied in \cite{Bo} and \cite{Vekua}. The solvability of the Dirichlet problem in  the partial case, when $\nu(z)=0$ and
the degeneracy of the ellipticity for the Beltrami equations
\begin{equation}\label{Beltrami1}
f_{\overline{z}}\, =\, \mu
(z)\cdot{f_z}
\end{equation}
 is controlled by the dilatation coefficient
\begin{equation}
 K_\mu(z)=K_{\mu,0}(z)={1+|\mu(z)|\over 1-|\mu(z)|}\notin\
L^{\infty},
\end{equation}
 is given in \cite{Dy}, \cite{GRSY2} and \cite{KPR1}.

\par

Recall that the problem on existence of homeomorphic solutions for
the equation (\ref{Beltrami1}) was resolved for the uniformly
elliptic case when $\Vert\mu\Vert_{\infty}<1$ long ago,  see e.g.
\cite{Ah$_1$}, \cite{Bo}, \cite{LV}. The existence problem for the
degenerate Beltrami equations (\ref{Beltrami1}) when $K_{\mu}
\notin\ L^{\infty}$
 is currently an active area of research, see e.g. the monographs \cite{GRSY2} and
\cite{MRSY} and the surveys \cite{GRSY1} and \cite{SY} and further
references therein. A series of criteria on the existence of
regular solutions for the Beltrami equation (\ref{Beltrami}) were
given in our recent papers \cite{BGR1}--\cite{BGR3}. There we
called a homeomorphism $f\in W^{1,1}_{\rm loc}(D)$ by a {\bf
regular solution} of (\ref{Beltrami}) if $f$ satisfies
(\ref{Beltrami}) a.e. in $D$ and $J_f(z)=|f_z|^2-|f_{\bar z}|^2\ne
0$ a.e. in $D.$
\medskip

\medskip

\cc
\section{Preliminaries}

To derive  criteria for existence of regular solutions for the
Dirichlet problem (\ref{Dirichlet}) in a Jordan domain $D\in {\Bbb C}$
we make use of the approximate procedure based on
the existence theorems for the case $K_{\mu , \nu}\in L^{\infty}$
given in \cite{Bo} and convergence theorems for the Beltrami
equations (\ref{Beltrami}) when $K_{\mu , \nu}\in L^{1}_{\rm loc}$
established in \cite{BGR2}. The
Schwarz formula
\begin{equation}\label{eqKPRS1.3}f(z)=i\,{\rm
Im}\,f(0)+\frac{1}{2\pi i}\int\limits_{|\zeta|=1}{\rm
Re}\,f(\zeta)\cdot\frac{\zeta+z}{\zeta-z}\frac{d\zeta}{\zeta}\,,\end{equation}
that allows to recover an analytic function $f$ in the unit disk
$\B$ by its real part $\varphi(\zeta)={\rm Re}\,f(\zeta)$ on the
boundary of ${\B}$ up to a purely imaginary additive constant
$c=i{\rm Im}\,f(0),$ see, e.g., Section 8, Chapter III, Part 3 in
\cite{HuCo}, as well as the Arzela--Askoli theorem combined with
moduli techniques  are also used.

\medskip

The following statement, that is a consequence of  Theorems 5.1 and
6.1 and the point 8.1 in \cite{Bo}, is basic for our further
considerations. See also Theorem VI.2.2 and the point VI.2.3 in
\cite{LV}, on the regularity of a $W^{1,1}_{\rm loc}$ solution to
the Beltrami equation (\ref{Beltrami1})     with the bounded
dilatation coefficient $K_\mu .$

\bpr \label{pr2} Let $D,$ $0\in D,$ be a Jordan domain in the complex plane ${\Bbb C}$ and $\varphi:\partial D\to{\Bbb R}$ be a
nonconstant continuous function. If $K_{\mu ,
\nu}\in L^{\infty},$ then the Dirichlet problem
(\ref{Dirichlet}) has
the unique regular solution $f$ normalized by ${\rm Im}
f(0)=0.$
This solution has the representation \begin{equation}
\label{eq4} f={\cal A}\circ g\circ {\cal R} \end{equation} where ${\cal R}:
D\to\B,$ ${\cal R}(0)=0,$  is a conformal mapping  and $g: \overline{\B}\to\overline{\B}$ stands for a
homeomorphic regular solution of the quasilinear equation
\begin{equation} \label{Beltrami3} g_{\overline{\zeta}}\, =\, \mu^*
(\zeta)\cdot{g_{\zeta}} +\nu^* ( \zeta )\cdot \frac{\overline{{\cal
A}^{\prime}(g(\zeta))}}{{\cal A}^{\prime}(g(\zeta))}\cdot \overline
{g_{ \zeta }} \end{equation} in $\B$  normalized by  $g(0)=0,$ $g(1)=1.$ Here $\mu*=\frac{{{\cal R}^{\prime}}}{\overline{{\cal
R}^{\prime}}} \cdot\mu\circ {\cal R}^{-1},$  $\nu*=\nu\circ {\cal
R}^{-1}$ and
\begin{equation}\label{eqKPRS1.33}{\cal A}(w)\colon =\frac{1}{2\pi i}\int\limits_{|\omega|=1}
\varphi({\cal
R}^{-1}(g^{-1}(\omega)))\cdot\frac{\omega+w}{\omega-w}\frac{d\omega}{\omega}\end{equation}
is an analytic function in the unit disk $\B$. \epr

\brem\label{rmk11} Let $\tilde\mu:\C\to\C$ coincide a.e in the
domain $D$ with \begin{equation} \label{Beltrami33} \frac{(g\circ
{\cal R})_{\overline{z}}}{(g\circ{\cal R})_z}\, =\,
\frac{g_{\overline\zeta}\circ{\cal R}\cdot \overline{{\cal
R}^{\prime}}}{g_{\zeta}\circ{\cal R}\cdot{{\cal R}^{\prime}}} =\,
\mu  + \nu\cdot\frac{\overline{{\cal R}^{\prime}}}{{\cal
R}^{\prime}} \cdot \frac{\overline{g_{{\zeta}}}}{g_{\zeta}}\circ
{\cal R} \cdot \frac{\overline{{\cal A}^{\prime}}}{{\cal
A}^{\prime}}\circ g\circ {\cal R} \end{equation} and equal to $0$
outside of $D$, see e.g. the formulas I.C(1) in \cite{Ah$_1$}.
Note that $K_{\tilde\mu}\le K_{\mu , \nu}$ a.e. in $D$ and there
is a regular solution $G:\overline\C\to\overline\C$ of the
equation $G_{\overline{z}}=\tilde\mu G_z$ such that $G(0)=0$,
$|G({\cal R}^{-1}(1))|=1,$  $G(\infty)=\infty$ and $G={\cal H}\circ
g\circ {\cal R}$ in $\overline{D}.$ Here ${\cal H}:\B\to G(D)$ is
a conformal mapping normalized by ${\cal H}(0)=0,$  ${\cal
H}^{\prime}(0)>0$. Thus, \begin{equation} \label{eq44} f={\cal
A}\circ h,
\end{equation}
\begin{equation}\label{eqKPRS33}{\cal A}(w)=\frac{1}{2\pi i}\int\limits_{|\omega|=1}
\varphi(h^{-1}(\omega))\cdot\frac{\omega+w}{\omega-w}\frac{d\omega}{\omega}\end{equation}
where \begin{equation} \label{eq444}
 h=g\circ{\cal R}={\cal H}^{-1}\circ G\end{equation} stands for a
homeomorphism $ h:\overline{D}\to\overline{\B}$, $h(0)=0$, which is
a regular solution in $D$ of the quasilinear equation \begin{equation}
\label{eq111.33} h_{\overline{z}}\, =\, \mu (z)\cdot h_z +\nu
(z)\cdot\frac{\overline{{\cal A}^{\prime}(h(z))}}{{\cal
A}^{\prime}(h(z))}\cdot \overline {h_z} \end{equation}

Denote such $f$, $g$, $\cal A$, $G$, $\cal H$ and $h$ by $f_{\mu ,
\nu , \varphi}$, $g_{\mu , \nu , \varphi}$, $\cal A_{\mu , \nu ,
\varphi}$ $G_{\mu , \nu , \varphi}$, ${\cal H}_{\mu , \nu ,
\varphi}$ and $h_{\mu , \nu , \varphi}$, respectively.\erem

Recall also that, given a family of paths $\Gamma $ in $\C ,$ a
Borel function $\rho:\C \to [0,\infty]$ is called {\bf admissible}
for $\Gamma ,$ abbr. $\rho \in adm\, \Gamma ,$ if \begin{equation}
\label{eq1.2v} \int\limits_{\gamma} \rho(z)\, |dz|\ \geq\ 1 \end{equation} for
each $\gamma\in\Gamma .$ The {\bf modulus} of $\Gamma$ is defined by
\begin{equation} \label{Beltramiv} M(\Gamma) =\inf\limits_{ \rho \in adm\, \Gamma}
\int\limits_{\C} \rho^2(z)\ dxdy\ . \end{equation}

\brem\label{rmk2.1} Note the following useful fact for a
quasiconformal mapping $f: D\to\C$, see e.g. V(6.6) in \cite{LV},
that \begin{equation} \label{eq2.14} M(f(\G))\ \le\ \int\limits_{\C} K(z)\cdot
\r^2(z)\ dxdy \end{equation} for every path family $\G$ in $D$ and for all
$\rho \in adm\, \Gamma$ where \begin{equation} \label{eq2.15} K(z)\ =\
\frac{|f_z|+|f_{\overline {z}}|}{|f_z|-|f_{\overline {z}}|} \end{equation} is
the (local) maximal dilatation of the mapping $f$ at a point $z\in
D.$

\erem

Given a domain $D$ and two sets $E$ and $F$ in ${\lC}$,
$\Delta(E,F,D)$ denotes the family of all paths $\g:[a,b] \to {\lC}$
which join $E$ and $F$ in $D$, i.e., $\g(a) \in E, \ \g(b) \in F$
and $\g(t) \in D$ for $a<t<b$. Recall that a {\bf ring domain}, or
shortly a {\bf ring} in $\lC$ is a domain $R$ whose complement
$\lC\setminus R$ consists of two connected components.
\medskip

Recall that, for points $z,\z\in\lC ,$ the {\bf spherical (chordal)
distance} $s(z,\z)$ between $z$ and $\z$ is given by \begin{equation}
\label{eq1.5a}
 s(z,\z )\  =\
\frac{|z-\z|}{(1+|z|^2)^{\frac{1}{2}}(1+|\z|^2)^{\frac{1}{2}}}\ \ \
{\mbox{if}}\ \ \  z\ \ne\ \infty\ne\z\ ,\end{equation}
$$ s(z,\infty )\ =\ \frac{1}{(1+|z|^2)^{\frac{1}{2}}}\ \ \
{\mbox{if}}\ \ \  z\ \ne\ \infty\ .$$ By $\d(A)$ we denote the
spherical diameter of a set $A\subset\C$, i.e. $\sup\limits_{z,\z\in
A}s(z,\z)$.

The following statement is a direct consequence of the known
estimate of the capacity of a ring formulated in terms of moduli,
see e.g. Lemma 2.16 in \cite{BGR2}.

\blem{} \label{lem4.1B} Let $f:D\to\C$ be a homeomorphism with $\d
(\lC\setminus f(D)) \ge \D > 0$ and let $z_0$ be a point in $D,$
$\z\in B(z_0,r_0),$ $r_0 < dist\, (z_0,\partial D).$ Then \begin{equation}
\label{eq4.33C} s(f(\z), f(z_0))\ \le\ \frac{32}{\D}\cdot
\hbox{exp}\left(-\frac{2\pi}{M(\D(fC,fC_0, fA))} \right) \end{equation} where
$C_0=\{z\in\C: |z-z_0|=r_0\}$, $C=\{z\in\C : |z-z_0|=|\z -z_0|\}$
and $A=\{z\in\C : |\z -z_0|<|z-z_0|< r_0\} .$\elem

\cc
\section{BMO, VMO and FMO functions}

Recall that a real-valued function $u$ in a domain $D$ in ${\Bbb C}$
is said to be of {\bf bounded mean oscillation} in $D$, abbr.
$u\in{\rm BMO}(D)$, if $u\in L_{\rm loc}^1(D)$ and
\begin{equation}\label{lasibm_2.2_1}\Vert u\Vert_{*}:=
\sup\limits_{B}{\frac{1}{|B|}}\int\limits_{B}|u(z)-u_{B}|\,dxdy<\infty\,,
\end{equation}
where the supremum is taken over all discs $B$ in $D$ and
$$u_{B}={\frac{1}{|B|}}\int\limits_{B}u(z)\,dxdy\,.$$ We write $u\in{\rm BMO}_{\rm loc}(D)$ if
$u\in{\rm BMO}(U)$ for every relatively compact subdomain $U$ of $D$
(we also write BMO or ${\rm BMO}_{\rm loc }$ if it is clear from the
context what $D$ is).

\medskip

The class BMO was introduced by John and Nirenberg (1961) in the
paper \cite{JN} and soon became an important concept in harmonic
analysis, partial differential equations and related areas, see e.g.
\cite{HKM} and \cite{RR}.

\medskip

A function $u$ in BMO is said to have {\bf vanishing mean
oscillation}, abbr. $u\in{\rm VMO}$, if the supremum in
(\ref{lasibm_2.2_1}) taken over all balls $B$ in $D$ with
$|B|<\varepsilon$ converges to $0$ as $\varepsilon\to0$. VMO has
been introduced by Sarason in \cite{Sarason}. There exists a number
of papers devoted to the study of partial differential equations
with coefficients of the class VMO.

\brem\label{rmk1} Note that $W^{\,1,2}\left({{D}}\right) \subset VMO
\left({{D}}\right),$ see e.g. \cite{BN}. \erem

Following \cite{IR}, we say that a function $u: D\to \R $ has {\bf
finite mean oscillation} at a point $z_0 \in {D} $ if
 \begin{equation} \label{eq2.4} \overline{\lim\limits_{\varepsilon\to 0}}\
 \ \
\dashint_{B( z_0 ,\varepsilon)}
|u(z)-\tilde{u}_{\varepsilon}(z_0)|\ dxdy\ <\ \infty \end{equation} where
$$ \tilde{u}_{\varepsilon}(z_0)=\dashint_{B( z_0 ,\varepsilon)}
u(z)\ dxdy $$ is the mean value of the function $u(z) $ over the
disk $B( z_0 ,\varepsilon)$ with small $\e>0.$ We also say that a
function $u : D\to \R $ is of {\bf finite mean oscillation } in
$D$, abbr. $u\in$ FMO(D) or simply $u\in$ {\bf FMO}, if
(\ref{eq2.4}) holds at every point $z_0 \in {D}.$

\brem\label{rmk2.33} Clearly BMO $\subset$ FMO. There exist examples
showing that FMO is not BMO$_{\rm loc},$ see e.g. \cite{GRSY2}. By
definition FMO$\ \subset L^1_{\rm loc}$ but FMO is not a subset of
$L^p_{\rm loc}$ for any $p>1$ in comparison with BMO$_{\rm loc}\subset
L^p_{\rm loc}$ for all $p\in [1,\infty)$. \erem

\bpr \label{pr2.1} If, for some collection of numbers
$u_{\varepsilon}\in {\R},\ \ \varepsilon \in (0,\varepsilon_0] $,
\begin{equation} \label{eq2.7} \overline{\lim\limits_{\varepsilon\to 0}}\ \ \
\dashint_{B( z_0 ,\varepsilon)} |u(z)-u_{\varepsilon}|\ dxdy\ <
\infty\, , \end{equation} then $u $ is of finite mean oscillation at
$z_0$. \epr

\bcor \label{cor2.1} If, for a point $z_0\in{D} ,$ \begin{equation}
\label{eq2.8} \overline{\lim\limits_{\varepsilon\to 0}}\ \
\dashint_{B( z_0 ,\varepsilon)} |u(z)|\ dxdy\ <\ \infty \ , \end{equation}
then $u $  has finite mean oscillation at $z_0.$ \ecor

\brem\label{rmk2.13a}  Note that the function
$u(z)=\log\frac{1}{|z|}$ belongs to BMO in the unit disk $\B$, see
e.g. \cite{RR}, p. 5, and hence also to FMO. However,
$\tilde{u}_{\e}(0)\to\infty$ as $\e\to 0,$ showing that the
condition (\ref{eq2.8}) is only sufficient but not necessary for a
function $u$ to be of finite mean oscillation at $z_0.$ \erem

Below we use the notation $
A(\varepsilon,\varepsilon_0)=\{z\in{\C}:\varepsilon<|z|<\varepsilon_0\}
\, .$

\blem{} \label{lem2.1}Let $u: D\rightarrow {\R}$ be a nonnegative
function with finite mean oscillation at $0\in {D}$ and let $u$
be integrable in $B(0,e^{-1})\subset D.$ Then \begin{equation} \label{eq2.9}
\int\limits_{A(\varepsilon, e^{-1})}\frac{u (z)\, dxdy}
{\left(|z| \log \frac{1}{|z|}\right)^2} \le\ C \cdot \log\log
\frac{1}{\varepsilon}\ \ \ \ \ \ \ \ \ \ \ \ \forall\  \e\in
(0,e^{-e}) \end{equation} \elem

For the proof of this lemma, see \cite{IR}.


\cc
\section{The main lemma}
The following lemma is the main tool for deriving   criteria on the
existence of regular solutions for the Dirichlet problem to the
Beltrami equations with two characteristics in a Jordan
domain in $\C$.

\blem{} \label{lem3.3A} Let $D$ be a Jordan domain in $\C$ with
$0\in D$ and let $\mu$ and $\nu : D\to\C$ be measurable functions
with $K_{\mu , \nu}\in L^1(D).$ Suppose that for every $z_0\in
\overline D $ there exist $\e_0=\e(z_0)>0$ and a family of
measurable functions $\p_{z_0,\e}:(0,\infty)\to(0,\infty),$
$\e\in(0,\e_0),$ such that \begin{equation} \label{eq3.5A} 0\ <\ I_{z_0}(\e)\
\colon =\ \int\limits_{\e}^{\e_0}\p_{z_0,\e}(t)\ dt\ <\ \infty\ ,
\end{equation}  and such that \begin{equation}\label{eq3.4o}
\int\limits_{\e<|z-z_0|<\e_0}\ K_{\mu ,
\nu}(z)\cdot\p^2_{z_0,\e}(|z-z_0|)\ dxdy\ =\ o(I^2_{z_0}(\e)) \end{equation}
as $\e\to 0.$ Then the  the Dirichlet problem
(\ref{Dirichlet}) has
a regular solution $f$
with ${\rm Im}
f(0)=0$ for each nonconstant continuous function $\f:\partial
D\to{\Bbb R}$. \elem

Here we assume that $\mu$ and $\nu$ are extended by zero outside of
the domain $D$.

\medskip

{\it Proof.} Setting  \begin{equation} \label{eq3.21p}\m_n(z)\ =\ \left
\{\begin{array}{rr} \m(z)\ , & \ \mbox{if}\  K_{\mu , \nu}(z)\le n,
\\ 0\ ,  &  \ \mbox{otherwise in}\ \C,
\end{array} \right. \end{equation} and \begin{equation} \label{eq3.21p}\n_n(z)\
=\ \left \{\begin{array}{rr} \n(z)\ , & \ \mbox{if}\ K_{\mu ,
\nu}(z)\le n,
\\ 0\ ,  &  \ \mbox{otherwise in}\ \C,
\end{array} \right. \end{equation} we have that $K_{\mu_n , \nu_n}(z)\le n$
in $\C$. Denote by $f_n$, ${\cal A}_n$, $G_n$, ${\cal H}_n$ and
$h_n$,  the functions $f_{\mu_n , \nu_n , \f}$,  ${\cal
A}_{{\mu }_n , {\nu }_n , \f}$ $G_{\mu_n , \nu_n , \f}$,
${\cal H}_{\mu_n , \nu_n , \f}$ and $h_{\mu_n , \nu_n ,
\f}$,, respectively, from Proposition \ref{pr2} and Remark
\ref{rmk11}.
\medskip

Let $\Gamma_{\varepsilon}$ be a family of all paths joining the
circles $C_{\varepsilon}=\{ z\in\C:|z-z_0|=\varepsilon\}$ and
$C_0=\{ z\in\C:|z-z_0|=\varepsilon_0\}$ in the ring
$A_{\varepsilon}=\{ z\in\C:\varepsilon< |z-z_0|<\varepsilon_0\}$.
Let also $\p^*$ be a Borel function such that $\p^*(t)=\p (t)$ for
a.e. $t\in (0,\infty )$. Such a function $\p^*$ exists by the Lusin
theorem, see e.g. \cite{Sa}, p. 69. Then the function
$$\rho_{\varepsilon}(z)=\left\{\begin{array}{rr}
\p^*(|z-z_0|)/I_{z_0}(\varepsilon), &  {\rm if } \ z\in A_{\varepsilon}, \\
0,  & {\rm  if} \ z\in {\C}\backslash A_{\varepsilon},
\end{array}\right.$$ is admissible for $\Gamma_{\varepsilon}$.
Hence by Remark \ref{rmk2.1} applied to $G_n$
$$M(G_n\Gamma_{\varepsilon})\leq\int\limits_{\varepsilon<|z-z_0|<\varepsilon_0} K_{\m ,\n}(z)\cdot
{\rho_{\varepsilon}}^2 (|z-z_0|)\ dx dy\,,$$ and, by the condition
(\ref{eq3.4o}), $M(G_n\Gamma_{\varepsilon})\to 0$ as $\varepsilon\to
0$ uniformly with respect to the parameter $n=1,2,\dots$.

\medskip

Thus, in view of the normalization $G_n(0)=0,$ $|G_n({\cal
R}^{-1}(1))|=1$, $G_n(\infty)=\infty$, the sequence $G_n$ is
equicontinuous in $\overline\C$ with respect to the spherical
distance by Lemma \ref{lem4.1B} with $\D = 1/{\sqrt{2}}$.
Consequently, by the Arzela--Ascoli theorem, see e.g. \cite{Du}, p.
267, and \cite{DS}, p. 382, it has a subsequence $G_{n_l}$ which
converges uniformly in $\overline\C$ with respect to the spherical
metric to a continuous mapping $G$ in $\overline\C$ with the
normalization $G(0)=0,$ $|G({\cal R}^{-1}(1))|=1$,
$G(\infty)=\infty$. Note that $G:\overline\C\to\overline\C$ is a
homeomorphism of the class $W^{1,1}_{\rm loc}(\C)$ by Corollary 3.8 in
\cite{BGR2}.

\medskip

Hence by the Rado theorem, see e.g. Theorem II.5.2 in \cite{Go},
${\cal H}_{n_l}\to {\cal H}$ as $l\to\infty$ uniformly in
$\overline\B$ where ${\cal H}:\overline\B\to G(\overline D)$ is the
conformal mapping of $\B$ onto $G(D)$ with the normalization ${\cal
H}(0)=0$ and ${\cal H}^{\prime}(0)>0$. Moreover, since the locally
uniform convergence $G_{n_l} \to G$ and ${\cal H}_{n_l} \to {\cal
H}$ of the sequences $G_{n_l}$ and ${\cal H}_{n_l}$ is equivalent to
their continuous convergence, i.e., $G_{n_l}(z_l) \to G(z_*)$ if
$z_l \to z_*$ and ${\cal H}_{n_l}(\zeta_l) \to {\cal H}(\zeta_*)$ if
$\zeta_l \to \zeta_*$, see [Du], p. 268, and since $G$ and ${\cal
H}$ are injective, it follows that $G_{n_l}^{-1} \to G^{-1}$ and
${\cal H}_{n_l}^{-1} \to {\cal H}^{-1}$ continuously, and hence
locally uniformly.

\medskip

Then we have that ${\cal A}_{n_l}\to{\cal A}$ locally uniformly in
$\B$ where
\begin{equation}\label{eqKPRS333}{\cal A}(w)=\frac{1}{2\pi i}\int\limits_{|\omega|=1}
\f(h^{-1}(\omega))\cdot\frac{\omega+w}{\omega-w}\frac{d\omega}{\omega}\
\end{equation} where $h:\overline D\to\overline\B$, $h(0)=0$, is a homeomorphism
$h={\cal H}^{-1}\circ G$. Note that ${\cal A}_{n_l}$ and $\cal A$
are not constant and hence ${\cal A}^{\prime}_{n_l}$ and ${\cal
A}^{\prime}$ have only isolated zeros. The collection of all such
zeros is countable. Thus, by Theorem 3.1 and Corollary 3.8 in
\cite{BGR2} $h_{n_l}\to h$ locally uniformly in $D$ and $h$ is a
homeomorphic $W^{1,1}_{\rm loc}$ solution in $D$ of the quasilinear
equation
\begin{equation} \label{eq111.333} h_{\overline{z}}\, =\, \mu
(z)\cdot{h_z} +\nu (z)\cdot\frac{\overline{{\cal
A}^{\prime}(h(z))}}{{\cal A}^{\prime}(h(z))}\cdot \overline {h_z}
\end{equation} Hence $f_{n_l}\to f$ where $f={\cal A}\circ h$ is a continuous
discrete open $W^{1,1}_{\rm loc}$ solution in $D$ of (\ref{Beltrami}).

\medskip

Next, note that ${\rm Re}\, {\cal A}_{n_l}\to{\rm Re}\, {\cal A}$
uniformly in $\overline\B$ by the maximum principle for harmonic
functions and ${\rm Re}\, {\cal A}=\f\circ h^{-1}$ on
$\partial\B$ and, consequently, ${\rm Re}\, f_{n_l}\to{\rm Re}\, f$
uniformly in $\overline\B$ and ${\rm Re}\, f=\f$ on $\partial
D$, i.e., $f$ is a continuous discrete open $W^{1,1}_{\rm loc}$ solution
of the Dirichlet problem (\ref{Dirichlet}) in $\B$ to the equation
(\ref{Beltrami}). It remains to show that $J_f(z)\ne 0$ a.e. in $\B$.

\medskip

By a change of variables which is permitted because $h_{n_l}$ and
$\tilde h_{n_l}=h_{n_l}^{-1}$ belong to the class $W^{1,2}_{\rm loc}$,
see e.g. Lemmas III.2.1 and III.3.2 and Theorems III.3.1 and III.6.1
in \cite{LV}, we obtain that for large enough $l$ \begin{equation}\label{eq4.2A}
\int\limits_{B} |\partial \tilde h_{n_l}|^2\ dudv\ \le
\int\limits_{\tilde h_{n_l}(B)} \frac{dxdy}{1-k_l(z)^2}\ \leq
\int\limits_{B^*} K_{\mu ,\nu }(z)\ dxdy\ <\ \infty\end{equation} where
$k_l(z)=|\mu_{n_l}(z)|+|\nu_{n_l}(z)|$ and $B^*$ and $B$ are
relatively compact domains in $D$ and $\tilde h(D)$, respectively,
such that $\tilde h(\bar{B}) \subset B^*$. The relation
(\ref{eq4.2A}) implies that the sequence $\tilde h_{n_l}$ is bounded
in W$^{1,2}(B)$, and hence $h^{-1} \in $ W$^{1,2}_{\rm loc},$ see e.g.
Lemma III.3.5 in \cite{Re} or Theorem 4.6.1 in \cite{EG}. The latter
condition brings in turn that $h$ has $(N^{-1})-$property, see e.g.
Theorem III.6.1 in \cite{LV}, and hence $J_h(z)\ne 0$ a.e., see
Theorem 1 in \cite{Po}. Thus, $f={\cal A}\circ h$ is a regular
solution of the Dirichlet problem (\ref{Dirichlet}) to the equation
(\ref{Beltrami}).

\bcor \label{cor4.3A} Let $D$ be a Jordan domain in $\C$ with $0\in
D$ and let $\mu$, $\nu : \B\to\C$ be measurable functions with
$K_{\mu , \nu}\in L^1(\B).$  Suppose that for every
$z_0\in\overline\B$ and some $\e_0>0$ \begin{equation} \label{eq3.4B}
\int\limits_{\e<|z-z_0|<\e_0} K_{\mu ,\nu}(z)\cdot\p^2(|z-z_0|)\
dxdy\ \le\ O\left( \int\limits_{\e}^{\e_0}\ \p(t)\ dt\right)\end{equation} as
$\e\to 0$, where $\p:(0,\infty)\to(0,\infty)$ is a measurable
function such that
\begin{equation} \label{eq3.5B}
\int\limits_{0}^{\e_0}\p(t)\
dt\ =\ \infty\ , \ \ \  0\ <\ \int\limits_{\e}^{\e_0}\p(t)\
dt
 <\ \infty\  \ \ \ \forall\ \e\in(0,\e_0) \ .
 \end{equation}
 Then the  the Dirichlet problem
(\ref{Dirichlet}) has
a regular solution $f$
with ${\rm Im}
f(0)=0$ for each nonconstant continuous function $\f:\partial
D\to{\Bbb R}$.
 \ecor

\cc
\section{Existence theorems}

Everywhere further we assume that the functions $\mu$ and $\nu :
D\to\C$ are extended by zero outside of the domain $D$.

\bth{} \label{th4.111a} Let $D$ be a Jordan domain in $\C$ with
$0\in D$ and let $\mu$ and $\nu :D\to\C$ be measurable functions
such that  $K_{\mu , \nu}(z)\,\ \le\ Q(z)\ \in\ \hbox{FMO}.$
Then the  the Dirichlet problem
(\ref{Dirichlet}) has
a regular solution $f$
with ${\rm Im}
f(0)=0$ for each nonconstant continuous function $\f:\partial
D\to{\Bbb R}$.  \eth

{\it Proof}. Lemma \ref{lem3.3A} yields this conclusion by
choosing \begin{equation} \label{eq3.4E} \p_{z_0,\e}(t)\  =\ \frac{1}{t\log
\frac{1}{t}}\ \ , \end{equation} see also Lemma \ref{lem2.1}.

\bcor \label{cor4.333b} In particular, if \begin{equation} \label{eq2.8a}
\overline{\lim\limits_{\varepsilon\to 0}}\ \ \ \dashint_{B( z_0
,\varepsilon)} \frac{1+|\nu (z)|}{1-|\nu (z)|}\ dxdy\ <\ \infty\ \ \
\ \ \ \ \ \ \forall\ z_0\in \overline D\ , \end{equation}
Then the  the Dirichlet problem
\begin{equation}\label{Dirichlet1}
\left\{\begin{array}{ccc}
f_{\overline{z}}\, =\, \nu (z)\cdot \overline {f_z},\,\,\, &z\in D, \\
\lim\limits_{z\to\zeta}{\rm Re}\,f(z)=\f(\zeta), &\forall\
\zeta\in\partial D,
\end{array}\right.
\end{equation}
 in a Jordan domain $D$, $0\in D,$ has
a regular solution $f$
with ${\rm Im}
f(0)=0$ for each nonconstant continuous function $\f:\partial
D\to{\Bbb R}$.
\ecor

Similarly, choosing in Lemma \ref{lem3.3A} the function
$\psi(t)=1/t$, we come to the following statement.

\medskip

\bth{} \label{thKPRS12b*} Let $D$ be a Jordan domain in $\C$ with
$0\in D$ and let $\mu$ and $\nu :D\to\C$ be measurable functions
such that  $K_{\mu , \nu}\in\ L^1_{\rm loc}(D).$ Suppose that
\begin{equation}\label{eqKPRS12c*}
\int\limits_{\varepsilon<|z-z_0|<\varepsilon_0}K_{\mu ,\nu}(z)\
\frac{dm(z)}{|z-z_0|^2}\ =\
o\left(\left[\log\frac{1}{\varepsilon}\right]^2\right)\qquad\forall\
z_0\in\overline{D}\end{equation} as $\varepsilon\to 0$ for some
$\varepsilon_0=\delta(z_0)$. Then the  the Dirichlet problem
(\ref{Dirichlet}) has a regular solution $f$ with ${\rm Im} f(0)=0$
for each nonconstant continuous function $\f:\partial D\to{\Bbb R}$.
\eth

\brem\label{rmKRRSa*} Choosing in Lemma \ref{lem3.3A} the function
$\psi(t)=1/(t\log{1/t})$ instead of $\psi(t)=1/t$, we are able to
replace (\ref{eqKPRS12c*}) by
\begin{equation}\label{eqKPRS12f*}
\int\limits_{\varepsilon<|z-z_0|<\varepsilon_0}\frac{K_{\mu
,\nu}(z)\ dm(z)} {\left(|z-z_0|\log{\frac{1}{|z-z_0|}}\right)^2}
=o\left(\left[\log\log\frac{1}{\varepsilon}\right]^2\right)\end{equation}
In general, we are able to give here the whole scale of the
corresponding conditions in $\log$ using functions $\psi(t)$ of the
form
$1/(t\log{1}/{t}\cdot\log\log{1}/{t}\cdot\ldots\cdot\log\ldots\log{1}/{t})$.
\erem

\bth{} \label{th3.2C}  Let $D$ be a Jordan domain in $\C$ with $0\in
D$ and let $\mu$, $\nu : D\to\B$ be measurable functions, $K_{\mu ,
\nu}\in L^1(D)$ and $k_{z_0}(r)$ be the mean value of $K_{\mu ,
\nu}(z)$ over the circle $|z-z_0|=r.$ Suppose that
 \begin{equation} \label{eq1}
 \int\limits_{0}^{\d(z_0)}\frac{dr}{rk_{z_0}( r)}\ =\ \infty
 \ \ \  \ \ \ \ \ \ \forall\ z_0\in \overline D\ .  \end{equation}
Then the  the Dirichlet problem
(\ref{Dirichlet}) has
a regular solution $f$
with ${\rm Im}
f(0)=0$ for each nonconstant continuous function $\f:\partial
D\to{\Bbb R}$. \eth

{\it Proof}. Theorem \ref{th3.2C} follows from Lemma \ref{lem3.3A}
by special choosing the functional parameter \begin{equation}
\label{eq3.21p}\p_{z_0,\e}(t)\ \equiv\ \p_{z_0}(t)\ \colon   =\
\left \{\begin{array}{rr} 1/[tk_{z_0}(t)]\ , & \ t\in (0,\e_0)\ ,
\\ 0\ ,  &  \ \mbox{otherwise}
\end{array} \right. \end{equation} where $\e_0=\d(z_0).$

\bcor \label{cor3.2V} In particular, the conclusion of Theorem
\ref{th3.2C} holds if \begin{equation} \label{eq3.21W} k_{z_0}(r)\ =\ O\left(
\log \frac{1}{r}\right)\ \ \ \mbox{as}\ \ \ r\to 0\ \ \  \ \ \ \ \ \
\forall\ z_0\in \overline D\ .  \end{equation} \ecor

In fact, it is clear that the condition (\ref{eq1}) implies the
whole scale of conditions in terms of $\log$ with using in the right
hand side in (\ref{eq3.21W}) functions  of the form
$\log{1}/{r}\cdot\log\log{1}/{r}\cdot\ldots\cdot\log\ldots\log{1}/{r}$.

\bigskip

In the theory of  mappings called quasiconformal in the mean,
conditions of the type
\begin{equation}\label{eq2} \int\limits_{{D}} \Phi
(Q(z))\ dxdy\  <\ \infty\end{equation} are standard for various
characteristics  of these mappings.
In this connection, in the paper \cite{RSY1}, see also the monograph
\cite{GRSY2}, it was established the equivalence of various integral
conditions on the function $\Phi$. We
give here the  conditions for $\Phi$ under which
(\ref{eq2}) implies (\ref{eq1}).\medskip

Further we use the following notion of the inverse function for
mo\-no\-to\-ne functions. Namely, for every non-decreasing function
$\Phi:[0,\infty ]\to [0,\infty ] ,$ the {\bf inverse function}
$\Phi^{-1}:[0,\infty ]\to [0,\infty ]$ can be well defined by
setting
\begin{equation}\label{eq5.5CC} \Phi^{-1}(\tau)\ =\
\inf\limits_{\Phi(t)\ge \tau}\ t\ .
\end{equation} As usual, here $\inf$ is equal to $\infty$ if the set of
$t\in[0,\infty ]$ such that $\Phi(t)\ge \tau$ is empty. Note that
the function $\Phi^{-1}$ is non-decreasing, too.

\brem\label{rmk3.333} It is evident immediately by the definition
that \begin{equation}\label{eq5.5CCC} \F^{-1}(\F(t))\ \le\ t\ \ \ \ \ \ \ \
\forall\ t\in[ 0,\infty ] \end{equation} with the equality in (\ref{eq5.5CCC})
except intervals of constancy of the function $\Phi$. \erem

Further, in (\ref{eq333Y}) and (\ref{eq333F}), we complete the
definition of integrals by $\infty$ if $\F(t)=\infty ,$
correspondingly, $H(t)=\infty ,$ for all $t\ge T\in[0,\infty) .$ The
integral in (\ref{eq333F}) is understood as the Lebesgue--Stieltjes
integral and the integrals (\ref{eq333Y}) and
(\ref{eq333B})--(\ref{eq333A}) as the ordinary Lebesgue integrals.

\bpr \label{pr4.1aB} Let $\F:[0,\infty ]\to [0,\infty ]$ be a
non-decreasing function and set \begin{equation}\label{eq333E} H(t)\ =\ \log
\F(t)\ .\end{equation}

Then the equality \begin{equation}\label{eq333Y} \int\limits_{\D}^{\infty}
H'(t)\ \frac{dt}{t}\ =\ \infty   \end{equation} implies the equality
\begin{equation}\label{eq333F} \int\limits_{\D}^{\infty} \frac{dH(t)}{t}\ =\
\infty  \end{equation} and (\ref{eq333F}) is equivalent to \begin{equation}\label{eq333B}
\int\limits_{\D}^{\infty}H(t)\ \frac{dt}{t^2}\ =\ \infty \end{equation} for
some $\D>0,$ and (\ref{eq333B}) is equivalent to every of the
equalities: \begin{equation}\label{eq333C}
\int\limits_{0}^{\d}H\left(\frac{1}{t}\right)\ {dt}\ =\ \infty \end{equation}
for some $\d>0,$ \begin{equation}\label{eq333D} \int\limits_{\D_*}^{\infty}
\frac{d\eta}{H^{-1}(\eta)}\ =\ \infty \end{equation} for some $\D_*>H(+0),$
\begin{equation}\label{eq333A} \int\limits_{\d_*}^{\infty}\ \frac{d\t}{\t
\F^{-1}(\t )}\ =\ \infty \end{equation} for some $\d_*>\F(+0).$
\medskip

Moreover, (\ref{eq333Y}) is equivalent  to (\ref{eq333F}) and hence
(\ref{eq333Y})--(\ref{eq333A})
 are equivalent each to other  if $\F$ is in addition absolutely continuous.
In particular, all the conditions (\ref{eq333Y})--(\ref{eq333A}) are
equivalent if $\F$ is convex and non--decreasing. \epr

Finally, we give the connection of the above conditions with the
condition of the type (\ref{eq1}).
\medskip

Recall that a function  $\psi :[0,\infty ]\to [0,\infty ]$ is called
{\bf convex} if $\psi (\lambda t_1 + (1-\lambda) t_2)\le\lambda\psi
(t_1)+ (1-\lambda)\psi (t_2)$ for all $t_1$ and $t_2\in[0,\infty ]$
and $\lambda\in [0,1]$.\medskip

\bpr \label{th5.555} Let $Q:\B\to [0,\infty ]$ be a measurable
function such that \begin{equation}\label{eq5.555} \int\limits_{\B} \F (Q(z))\
dxdy\  <\ \infty\end{equation} where $\F:[0,\infty ]\to [0,\infty ]$ is a
non-decreasing convex function such that \begin{equation}\label{eq3.333a}
\int\limits_{\d}^{\infty}\ \frac{d\t}{\t \F^{-1}(\t )}\ =\ \infty
\end{equation} for some $\d\ > \F(0).$ Then \begin{equation}\label{eq3.333A}
\int\limits_{0}^{1}\ \frac{dr}{rq(r)}\ =\ \infty \end{equation} where $q(r)$
is the average of the function $Q(z)$ over the circle $|z|=r$. \epr

Finally, combining Propositions \ref{pr4.1aB} and \ref{th5.555} we
obtain the following conclusion.

\bcor \label{cor555} If $\F:[0,\infty ]\to [0,\infty ]$ is a
non-decreasing convex function and $Q$ satisfies the condition
(\ref{eq5.555}), then every of the conditions
(\ref{eq333Y})--(\ref{eq333A}) implies (\ref{eq3.333A}). \ecor

Immediately on the basis of Theorem \ref{th3.2C} and Corollary
\ref{cor555}, we obtain the next significant result.

\bth{} \label{th4.111a} Let  $D$ be a Jordan domain in $\C$ with
$0\in D$ and let $\mu$ and $\nu : D\to\C$ be measurable functions
such that \begin{equation}\label{eq4.2} \int\limits_{{D}} \Phi (K_{\m,\n}(z))\
dxdy\  <\ \infty \end{equation} where $\F:[0,\infty ]\to [0,\infty ]$ is a
non-decreasing convex function satisfying at least one of the
conditions (\ref{eq333Y})--(\ref{eq333A}). Then the  the Dirichlet problem
(\ref{Dirichlet}) has
a regular solution $f$
with ${\rm Im}
f(0)=0$ for each nonconstant continuous function $\f:\partial
D\to{\Bbb R}$.  \eth

On the same basis, we obtain the following consequence.

\bcor \label{cor333} In particular, the conclusion of Theorem
\ref{th4.111a} holds if
\begin{equation}\label{eqp.KR4.1c}\int\limits_{D\cap
U_{z_0}}e^{\alpha(z_0) K_{\mu , \nu}(z)}\,dxdy\ <\infty
\qquad\forall\ z_0\in \overline{D}\end{equation} for some
$\alpha(z_0)>0$
 and a neighborhood $U_{z_0}$ of the point $z_0$.
\ecor

\brem\label{rmk5.1} By the Stoilow theorem, see
e.g. \cite{Sto}, every regular solution $f$ to the Dirichlet problem
\begin{equation}\label{Dirichlet3}
\left\{\begin{array}{ccc}
f_{\overline{z}}\, =\, \mu (z)\cdot {f_z},\,\,\, &z\in D, \\
\lim\limits_{z\to\zeta}{\rm Re}\,f(z)=\f(\zeta), &\forall\
\zeta\in\partial D,
\end{array}\right.
\end{equation}
has the representation $f=h\circ g$ where $g:D\to\B$ stands for a homeomorphic $W^{1,1}_{\rm loc}$
solution to the Beltrami equation $g_{\overline{z}}\, =\, \mu (z)\cdot {g_z},$ and $h:\B\to{\Bbb C}$
is analytic. By Theorem 5.50 from  \cite{RSY1}
the conditions (\ref{eq333Y})--(\ref{eq333A}) are not only
sufficient but also necessary to have a homeomorphic $W^{1,1}_{\rm loc}$ solution
for all such Beltrami equations with the integral constraint
\begin{equation} \int\limits_{{D}} \Phi (K_{\m}(z))\
dxdy\  <\ \infty. \end{equation}

Note also that in the above theorem we may assume that the functions
$\F_{z_0}(t)$ and $\F(t)$ are not convex and non--decreasing on the
whole segment $[0,\infty]$ but only on a segment $[T,\infty]$ for
some $T\in(1,\infty)$. Indeed, every function
$\F:[0,\infty]\to[0,\infty]$ which is convex and non-decreasing on a
segment $[T,\infty]$, $T\in(0,\infty)$, can be replaced by a
non-decreasing convex function $\F_T:[0,\infty]\to[0,\infty]$ in the
following way. We set $\F_T(t)\equiv 0$ for all $t\in [0,T]$,
$\F(t)=\f(t)$, $t\in[T,T_*]$, and $\F_T\equiv \F(t)$,
$t\in[T_*,\infty]$, where $\t=\f(t)$ is the line passing through the
point $(0,T)$ and tangent to the graph of the function $\t=\F(t)$ at
a point $(T_*,\F(T_*))$, $T_*\ge T$. For such a function we have by
the construction that $\F_T(t)\le \F(t)$ for all $t\in[1,\infty]$
and $\F_T(t)=\F(t)$ for all $t\ge T_*$. \erem

The equation of the form \begin{equation} \label{eq6.1} f_{\overline{z}}\, =\,
\l (z)\ {\rm Re}\, f_z \end{equation} with $|\l(z)|<1$ a.e. is called a {\bf
reduced Beltrami equation}, considered e.g. in \cite{Bo$_3$} and
\cite{Vo}, though the term is not introduced there. The equation
(\ref{eq6.1}) can be written as the equation (\ref{Beltrami}) with \begin{equation}
\label{eq6.2} \m(z)\ =\ \n (z)\ =\ \frac{\l(z)}{2} \end{equation} and then
\begin{equation} \label{eq6.3} K_{\mu , \nu}(z)\ =\ K_{\l}(z)\ \colon =\
\frac{1+|\l (z)|}{1-|\l (z)|}\ .\end{equation} Thus, we obtain from Theorem
\ref{th4.111a} the following consequence for the reduced Beltrami
equations (\ref{eq6.1}).

\bth{} \label{th4.111aR} Let $D$ be a Jordan domain in $\C$ with
$0\in D$ and let $\l : D\to\C$ be a measurable function
 such that \begin{equation}\label{eq2.8aR}
\int\limits_{{D}} \Phi (K_{\l}(z))\ dxdy\ <\ \infty \end{equation} where
$\F:[0,\infty ]\to [0,\infty ]$ is a non-decreasing convex function
satisfying at least one of the conditions
(\ref{eq333Y})--(\ref{eq333A}).
Then the  the Dirichlet problem
\begin{equation}\label{Dirichlet2}
\left\{\begin{array}{ccc}
f_{\overline{z}}\, =\,
\l (z)\ {\rm Re}\, f_z ,\,\,\, &z\in D, \\
\lim\limits_{z\to\zeta}{\rm Re}\,f(z)=\f(\zeta), &\forall\
\zeta\in\partial D,
\end{array}\right.
\end{equation}
 in a Jordan domain $D$, $0\in D,$ has
a regular solution $f$
with ${\rm Im}
f(0)=0$ for each nonconstant continuous function $\f:\partial
D\to{\Bbb R}.$
\eth

Finally, on the basis of Corollary \ref{cor333}, we obtain the
following consequence.

\bcor \label{cor111} In particular, the conclusion of Theorem
\ref{th4.111aR} holds if
\begin{equation}\label{eqp.KR4.1c}\int\limits_{D\cap
U_{z_0}}e^{\alpha(z_0) K_{\lambda}(z)}\,dxdy\ <\infty \qquad\forall\
z_0\in \overline{D}\end{equation} for some $\alpha(z_0)>0$
 and a neighborhood $U_{z_0}$ of the point $z_0$.
\ecor

\brem\label{rmk6.1} Remarks \ref{rmk5.1} are valid for reduced
Beltrami equations. Moreover, the above results remain true for the
case in (\ref{Beltrami}) when \begin{equation} \label{eq6.14} \n(z)\ =\ \m (z)\
e^{i\theta (z)} \end{equation} with an arbitrary measurable function
$\theta(z): D\to\R$ and, in particular, for the equations of the
form \begin{equation} \label{eq6.1A} f_{\overline{z}}\, =\, \l (z)\ {\rm Im}\,
f_z \end{equation} with a measurable coefficient $\l : D\to\C$, $|\l(z)|<1$
a.e., see e.g. \cite{Bo$_3$}. \erem

Our approach makes possible, under the certain modification, to
obtain criteria on the existence of pseudoregular and multi-valued
solutions in finitely connected domains that will be published
elsewhere.

\noindent Bogdan Bojarski, Institute of Mathematics of Polish
Academy of Sciences,\\ ul. Sniadeckich 8, P.O. Box 21, 00--956
Warsaw, POLAND\\ Email: {\tt
B.Bojarski@impan.gov.pl}\\

\noindent Vladimir Gutlyanski, Vladimir Ryazanov, Inst. Appl. Math.
Mech., NASU,\\ 74 Roze Luxemburg str., 83114, Donetsk, UKRAINE,\\
Email: {\tt vladimirgut@mail.ru}, {\tt vl.ryazanov1@gmail.com}


\medskip
\begin{thebibliography}{99}
\small


\bibitem{Ah$_1$}
{\sc Ahlfors L.V.}, {\sl Lectures on Quasiconformal Mappings}, D.
Van Nostrand Company, Inc., Princeton etc., 1966.

\bibitem{Bo}{\sc Bojarski B.},
{\sl Generalized solutions of a system of differential equations of
the first order of the elliptic type with discontinuous
coefficients}, Mat. Sb. 43(85) (1957), no. 4, 451--503 [in Russian];
transl. in Rep. Univ. Jyv\"askyl\"a Dept. Math. Stat. 118 (2009),
1--64.


\bibitem{Bo$_3$}{\sc Bojarski B.},
{\sl Primary solutions of general Beltrami equations}, Ann. Acad.
Sci. Fenn. Math. 32 (2007), no. 2, 549--557.

\bibitem{BGR1}{\sc Bojarski B., Gutlyanskii V. and Ryazanov V.},
{\sl General Beltrami equations and BMO}, Ukr. Mat. Visn.
\textbf{5}(3)  (2008), 305-326; transl. in Ukrainian Math. Bull.
\textbf{5}(3) (2008), 305--326.

\bibitem{BGR2}
{\sc Bojarski~B., Gutlyanskii~V. and Ryazanov~V.}, {\sl On Beltrami
equations with two characteristics}, Comp. Var. Ell. Equ. 54 (2009),
933--950.

\bibitem{BGR3}
B. Bojarski, V. Gutlyanskii and  V. Ryazanov, \textit{On Integral
Conditions for the General Beltrami Equations.} Comp. Anal. Oper.
Theory, \textbf{5} (3) (2011), 835--845.

\bibitem{BN}
 H. Brezis, L. Nirenberg, \textit{Degree theory and BMO. I. Compact manifolds
without boundaries.} Selecta Math. (N.S.) \textbf{1}(2) (1995),
197--263.

\bibitem{Du}
{\sc Dugundji J.}, {\sl Topology}, Allyn and Bacon, Inc., Boston,
1966.

\bibitem{DS}
{\sc Dunford N. and Schwartz J.T.}, {\sl Linear Operators, Part I:
General Theory}, Interscience Publishers Inc., New York, London,
1957.

\bibitem{Dy}
Yu. Dybov, \textit{On regular solutions of the Dirichlet problem for
the Beltrami equations.} Comp. Var. Ell. Equ. \textbf{55} (12)
(2010), 1099-�1116.

\bibitem{EG}
{\sc Evans L.C. and Gapiery R.F.}, {\sl Measure Theory and Fine
Properties of Functions}. CRC Press, Boca Raton, FL (1992).

\bibitem{Go}
G.M. Goluzin, \textit{Geometric Theory of Functions of a Complex
Variable} [in Russian], 2nd ed., Izdat. Akad. Nauk SSSR, 1966; Engl.
transl. in Amer. Math. Soc. Transl., v. 26, 1969, pp. 136--140,
165--170, 375�-379, 599, 643�-644.


\bibitem{GRSY1}
V. Gutlyanskii, V. Ryazanov,  U. Srebro, E. Yakubov, \textit{On
recent advances in the degenerate Beltrami equations.} Ukr. Mat.
Visn. \textbf{7}(4) (2010), 467--515; transl. in J. Math. Sci.,
\textbf{175}(4) (2011), 413--449.


\bibitem{GRSY2}
V. Gutlyanskii, V. Ryazanov, U. Srebro, E. Yakubov, ``The Beltrami
Equation: A Geometric Approach'', Developments in Mathematics,
\textbf{26}, Springer, New York, 2012.

\bibitem{HKM}
J. Heinonen, T. Kilpelainen, O. Martio, \textit{Nonlinear Potential
Theory of Degenerate Elliptic Equations.} Clarendon Press, Oxford
Univ. Press, 1993.

\bibitem{HuCo}
 A. Hurwitz, R. Courant, \textit{The Function Theory.} Nauka, Moscow, 1968 [in
Russian].

\bibitem{IR} {\sc Ignat'ev A. and Ryazanov V.}, {Finite mean
oscillation in the mapping theory}, Ukr. Mat. Visn. 2 (2005), no. 3,
395--417 [in Russian]; transl. in Ukrainian Math. Bull. 2 (2005),
no. 3, 403--424 .

\bibitem{JN}
{\sc John F. and Nirenberg L.}, {\sl On functions of bounded mean
oscillation}, Comm. Pure Appl. Math. 14 (1961), 415--426.

\bibitem{KK}{\sc Krushkal' S.L. and K\"uhnau R.},
{\sl Quasiconformal mappings: new methods and applications}, Nauka,
Novosibirsk, 1984 [in Russian].

\bibitem{KPR1}
D. Kovtonyuk, I. Petkov, V. Ryazanov, \textit{On the Dirichlet
problem for the Beltrami equations in finitely connected domains.}
Ukr. Mat. Zh. \textbf{64} (2012), 932--944.

\bibitem{LV}
{\sc Lehto O. and Virtanen K.}, {\sl Quasiconformal Mappings in the
Plane}, Springer, New York etc., 1973.


\bibitem{MRSY}
O. Martio, V. Ryazanov, U. Srebro, E. Yakubov, \textit{Moduli in
Modern Mapping Theory}, Springer Monographs in Mathematics,
Springer, New York etc., 2009.

\bibitem{Po} {\sc Ponomarev S.P.}, {\sl The $N^{-1}$--property of mappings,
and Lusin's (N) condition}, Mat. Zametki 58 (1995), 411-418; transl.
in Math. Notes 58 (1995), 960--965.

\bibitem{RR}{\sc Reimann H.M. and Rychener T.},
{\sl Funktionen Beschrankter Mittlerer Oscillation}, Springer,
Berlin etc., 1975.

\bibitem{Re}{\sc Reshetnyak Yu.G.},
{\sl Space Mappings with Bounded Distortion}, Transl. of Math.
Monographs 73, AMS, 1989.

\bibitem{RSY1}
V. Ryazanov, U. Srebro, E. Yakubov, \textit{Integral conditions in
the theory of the Beltrami equations.} Complex Var. Elliptic Equ. --
Published electronically May 2011. -- DOI:
10.1080/17476933.2010.534790 (expected in print, \textbf{57}, no.
12, 2012).

\bibitem{Sa}{\sc S. Saks},
{\sl Theory of the Integral}, Dover Publ. Inc., New York, 1964.

\bibitem{Sarason}
D. Sarason, \textit{Functions of vanishing mean oscillation.} Trans.
Amer. Math. Soc. \textbf{207} (1975), 391--405.

\bibitem{SY}
{\sc Srebro U. and Yakubov E.}, {\sl The Beltrami equation},
Handbook in Complex Analysis: Geometric function theory, Vol. 2,
555-597, Elseiver B. V., 2005.

\bibitem{Sto}
S. Stoilow, \textit{Lecons sur les Principes Topologue de le Theorie
des Fonctions Analytique.} Gauthier-Villars, 1938. Riemann,
Gauthier-Villars, Paris, 1956 [in French].

\bibitem{Vekua}
I.N. Vekua, \textit{Generalized analytic functions.} Pergamon Press,
London, 1962.

\bibitem{Vo}{\sc Volkovyskii, L. I.},
{\sl Quasiconformal mappings}, L'vov Univ. Press, L'vov, 1954 [in
Russian].

\end{thebibliography}
\end{document}